# Dynamic Cooperative Vehicle Platoon Control Considering Longitudinal and Lane-changing Dynamics

Kangning Hou, Fangfang Zheng, Xiaobo Liu, Zhichen Fan

*Abstract*—This paper presents a distributed cascade Proportional Integral Derivate (DCPID) control algorithm for the connected and automated vehicle (CAV) platoon considering the heterogeneity of CAVs in terms of the inertial lag. Furthermore, a real-time dynamic cooperative lane-changing model for CAVs, which can seamlessly combine the DCPID algorithm and the improved sine function is developed. The DCPID algorithm determines the appropriate longitudinal acceleration and speed of the lane-changing vehicle considering the speed fluctuations of the front vehicle on the target lane (TFV). In the meantime, the sine function plans a reference trajectory which is further updated in real time using the model predictive control (MPC) to avoid potential collisions until lane-changing is completed. Both the local and the asymptotic stability conditions of the DCPID algorithm are mathematically derived, and the sensitivity of the DCPID control parameters under different states is analyzed. Simulation experiments are conducted to assess the performance of the proposed model and the results indicate that the DCPID algorithm can provide robust control for tracking and adjusting the desired spacing and velocity for all 400 scenarios, even in the relatively extreme initial state. Besides, the proposed dynamic cooperative lane-changing model can guarantee an effective and safe lane-changing with different speeds and even in emergency situations (such as the sudden deceleration of the TFV).

## I. INTRODUCTION

WITH the rapid development of connected and automated vehicle (CAV) technologies, the transition from the traditional human-driven traffic to the CAV traffic is undergoing. CAVs can make full use of connectivity, especially fast V2V and V2I communication, to improve traffic capacity, safety and stability [1]–[3]. Among the applications of CAV technologies, the vehicle platooning has gained a lot of attention in recent years. Vehicle platooning is a coordinated movement mechanism where vehicles travel with small headways without any mechanical linkage [4]. The main benefits of vehicle platooning include energy consumption savings due to less aerodynamic drag [5], increased road capacity and enhanced traffic safety by providing lower reaction time [6], [7]. Therefore, it is of great significance to investigate the automatic control of vehicles for platooning. [1]

Basically, the vehicle platooning is mainly composed of two aspects: the longitudinal tracking control and the lateral lane-changing control. In the literature, considerable research has been carried out to develop longitudinal control models, among which the cooperative adaptive cruise control (CACC) model developed by the PATH laboratory, the model predictive control (MPC) models and the linear control models have gained substantial attention in both traffic and control engineering. The CACC technologies enable a vehicle to follow the preceding vehicle smoothly while maintaining a safety distance based on information obtained from the preceding vehicles, thus to improve the traffic efficiency and reduce oscillations [8]–[10]. Milanés *et al.* [8] first develop the CACC system by introducing V2V communications to commercially available ACC systems. The field road tests show that the CACC system could improve highway capacity and string stability. They further compare the CACC with the Intelligent Drive Model (IDM) and a commercial ACC through field experiments, and the test results indicate that the CACC performs best in terms of providing smooth and stable car following response [11]. Later, Liu *et al.* [12] extend the CACC model to consider both car following and lane changing behavior of CAVs in mixed traffic. However, the major drawback of the CACC model is that the control performance is highly dependent on scenarios and parameter settings.

The MPC models aim at optimizing the control decisions (e.g., acceleration/deceleration) of the following vehicles in the platoon for a certain prediction horizon based on the vehicles' prevailing state information [13]. Gong *et al.* [14] model the platoon as a multi-agent interconnected dynamic system and

Manuscript received July 16, 2021; This work was funded by National Science Foundation of China (NSFC) under project codes 61673321 and 52072315, Science & Technology Department of Sichuan Province under project codes 2019JDTD0002 and 2020JDJQ0034, and Chengdu Science and technology Bureau under project codes 2019-YF05-02657-SN. (Corresponding author: Fangfang Zheng)

K. Hou and F. Zheng are with School of Transportation and Logistics, Southwest Jiaotong University, Western Hi-tech Zone Chengdu, Sichuan 611756, P.R. China, and also with National Engineering Laboratory of Integrated Transportation Big Data Application Technology and National United Engineering Laboratory of Integrated and Intelligent Transportation, Southwest Jiaotong University, Western Hi-tech Zone Chengdu, Sichuan 611756, P.R. China. (e-mail: kangning.hou@my.swjtu.edu.cn; fzheng@swjtu.cn).
X. Liu and Z. Fan are with School of Transportation and Logistics, Southwest Jiaotong University, Western Hi-tech Zone Chengdu, Sichuan 611756, P.R. China. (e-mail: xiaobo.liu@swjtu.cn; fzc@swjtu.edu.cn).

propose a one-step MPC control algorithm based on constrained optimization and distributed computation. Gong and Du [15] further extend their model to the cooperative MPC-based platoon control in mixed traffic flow with CAVs and human-driven vehicles (HDVs) where system optimizers are developed to consider multiple objectives including traffic efficiency and driving/riding comfort. Liu *et al.* [16] formulate the platoon as a dynamically decoupled system and propose a two-step distributed MPC algorithm to solve the dynamic formation problem. The aforementioned platoon control strategies require the optimal control problem be solved instantaneously at each sampling time instant without considering the control delay. To address the problem of control delay, Wang *et al.* [13] develop two real-time deployable approaches based on the MPC which can provide efficient cooperative control for all following vehicles in the platoon to damp traffic oscillations. One significant advantage of the MPC model is that it can deal with multiple criteria and constraints on vehicles' state and control variables. However, the MPC approach, especially the centralized MPC which involves vehicles in a platoon, requires large amount of computation time and the stability of the controlled platoon through mathematical proofs is sometimes quite challenging [17].

The linear control models usually deal with designing feedback control laws (feedback gain matrices) for the state variables of the platoon system under different information flow topologies. The control values of the following vehicles are then calculated ensuring the stability of the closed-loop platoon system [18], [19]. Guo and Yue [18] propose a linear controller by acquiring the information of the preceding vehicle and the leading vehicle of the platoon. Their approach can stabilize the platoon with robustness for a given level of disturbance attenuation; while Ghasemi *et al.* [20] design a decentralized linear controller using the information of the preceding vehicle and the following vehicle as inputs. Zheng *et al.* [21] summarize six kinds of commonly used information flow topologies and establish linear controller gains for vehicle platoons. They further investigate the scalability and the stability margin of vehicle platoons under the undirected topology and the bidirectional topology [19]. However, there always exists oscillation around the steady-state (with zero spacing error) for the linear controllers discussed above.

Besides the aforementioned linear controllers, the PID controller is also widely used for the longitudinal tracking control and the lateral lane-changing control to attain a steady-state of the vehicle platoon [22]–[24]. These PIDs are generally centralized controllers, which need to collect all vehicle information and perform large-scale optimization calculations. To reduce computation costs, it is preferable to consider the vehicle platoon as a multiagent system (MAS) where the distributed framework can be applied to design controllers. In this case, each distributed controller only performs local optimal control for one vehicle, e.g., see [25], where a distributed robust PID control approach is proposed to ensure good leader-tracking performances for a platoon of connected autonomous vehicles. It is worth noting that the PID controllers developed in [19]-[22] are all single-layer PID controllers which cannot deal with problems involving constraints and multi-objective optimization. To deal with these problems, Lui *et al.* [26] proposed a fully distributed PID control strategy which can be used to solve both the leader-tracking and the containment control problems without requiring the fulfillment of additional constraints on the control input matrix.

For the lane-changing trajectory planning, most existing research applies the curve interpolation method or the artificial potential field method. The curve interpolation method is to select an appropriate curve to fit the path, such as the polynomial curve, Bessel curve, B-spline curve, Sine curve, etc. The time-dependent polynomial curve trajectory planning method [27]–[29] is to design polynomial functions of the longitudinal and the lateral positions with respect to time, where a large number of parameters need to be calibrated according to the state assumptions at the beginning and the end of the lane-changing process. These assumptions (e.g., the acceleration at the beginning and the end of lane-changing is 0) are often inconsistent with reality. Yang *et al.* [30] propose a time-independent polynomial curve method to fix unreasonable assumptions. Nevertheless, the model parameters still lack physical significance and are difficult to compute to achieve a deterministic motion state [31]. Maekawa *et al.* [32] apply the cubic B-spline curve to generate the lane-changing trajectory. Their approach couldn't guarantee the optimal solution because too much attention is paid to the continuity of the curvature, while the maximum lateral acceleration generated by the lane-changing trajectory is difficult to control. Chen *et al.* [33] develop a quadratic Bessel curve-based lane-changing trajectory planning algorithm considering the safe distance of lane-changing and the riding comfort. This method could provide flexibility in selecting control points and modifying the lane-changing distance once determined. Recently, Lee et al. [5] develop a vehicle trajectory generation approach for truck platooning by combing Kalman filter-based vehicle state estimation, front-rear trajectory buffering and the 3$^{rd}$ polynomial curve fitting method. The simulation results show that their proposed approach could maintain string stability for trucking platooning in different conditions.

The artificial potential field method plans the vehicle path by assigning different potential functions to different types of obstacles and road structures [34]. Rasekhipour *et al.* develop a potential field-based model predictive path-planning controller by combining the artificial potential field method with vehicle dynamics. Their method can guarantee vehicle stability during the path-tracking process. Huang *et al.* [35] use the velocity information in the path planning and regulate the velocity within the predicted horizon of the current state to satisfy a feasible local trajectory. However, the artificial potential field method usually has the problems of target unreachability and falling into local minimum. In addition, the performance of this method becomes worse with the increase of the speed resulting in poor robustness [36]. By far, the majority of the lane-changing models proposed in the literature are static. Though Yang *et al.* [30] propose a dynamic lane-changing trajectory planning model, the model does not consider the cooperative maneuverability of vehicles during lane-changing, which could be inconsistent with the real environment.

In this paper, we propose a dynamic cooperative vehicle platooning control approach which is composed of a distributed cascade PID (DCPID)-based longitudinal control model and a cooperative lane-changing model. The main contributions of this paper are as follows.
1. We propose a distributed double-layer DCPID control algorithm by combining an inner loop controlling the velocity error and an outer loop controlling the spacing error with objectives that both errors are zero. Both local and asymptotic stability conditions of the CAV platooning system is mathematically derived and the parameter sensitivity analysis is provided to guarantee good performance of the system.
2. We develop a dynamic cooperative lane-changing model which consists of the DCPID control algorithm and an improved sine curve trajectory planning method. The first component can ensure safe lane-changing by instructing the speed of the lane-changing vehicle. The second component can plan a continuous path during which the passenger discomfort caused by lateral acceleration is taken into account. Furthermore, the feasible domain of the planned acceleration of the improved sine function is determined using the yaw-rate of the vehicle, and the parameters of the improved sine function are explicable with physical significance.
3. The proposed model explicitly considers the heterogeneous behavior of the CAVs in terms of the inertial lag. The uncertainty of lane-changing caused by the speed change of the front vehicle on the target lane, and the cooperative maneuver between the vehicle and the platoon on the target lane are both considered.

The rest of this paper is organized as follows: Section II presents the detailed description of the proposed dynamic cooperative vehicle platoon control model, and Section III provides the detailed mathematical analysis of the stability conditions for the DCPID algorithm. In Section IV, the parameter sensitivity analysis of the DCPID is discussed. Section V gives the numerical examples both for the longitudinal control and the lane-changing trajectory planning and the simulation results are analyzed and discussed. Section VI concludes the paper, and some discussion and the prospect of future work are provided. Throughout the rest of the paper, all vehicles are referred to as CAVs.

## II. DYNAMIC COOPERATIVE VEHICLE PLATOON CONTROL MODEL

### A. Description of the cooperative control framework

We consider a vehicle platoon traveling on a two-lane freeway as shown in Fig. 1. A vehicle on the adjacent lane plans to change lanes and merge into the platoon. *SV* denotes the vehicle which desires to change lanes and join the platoon on the target lane; *TFV* and *TRV* represent the vehicles in front and in rear of the desired position of the *SV* on the target lane, respectively. Both the *TFV* and the *TRV* are the following vehicles in the original platoon and the *LV* is the leading vehicle of the platoon.

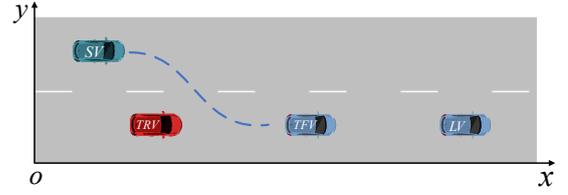

Fig.1. Illustration of vehicle lane-change

Fig. 2 illustrates the real-time dynamic cooperative vehicle platoon control framework which includes the longitudinal control model and the cooperative lane-changing model. The longitudinal control model is composed of a cascade PID control algorithm and a kinematic model, ensuring safe and efficient operation of the vehicle platoon. In case any vehicle in the platoon is disturbed, the longitudinal control model can respond quickly to make the platoon return to the steady state. The cooperative lane-changing model consists of three parts: lane-changing decision-making, dynamic trajectory generation and trajectory tracking.

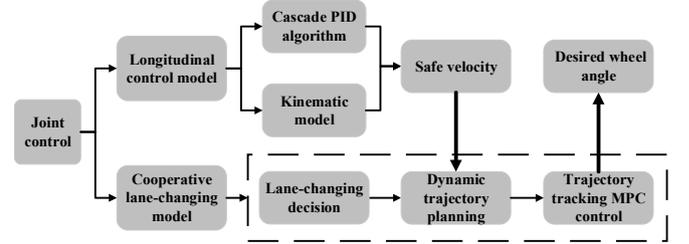

Fig. 2. The framework for real-time dynamic cooperative control

### B. Distributed Cascade PID algorithm for longitudinal control

Considering a platoon of *N* vehicles, we define any two adjacent vehicles on the same lane as a subsystem. Let $P_{i,i-1}$ denote the subsystem consisting of the preceding vehicle $i-1$ and the following vehicle $i$. The distributed cascade PID algorithm controls the operation of the following vehicle of each subsystem. Let $x_i(t), v_i(t), a_i(t)$ denote the position, the velocity and the acceleration of vehicle $i$ ($i \in \{1,2, ,N\}$) at time $t$, respectively. The longitudinal dynamics of the following vehicle in the subsystem is described by the linearized third-order state space equations as:

$$\begin{aligned}\dot{x}_i(t) &= v_i(t) \\ \dot{v}_i(t) &= a_i(t) \\ \dot{a}_i(t) &= -\frac{1}{\tau_i} a_i(t) + \frac{1}{\tau_i} u_i(t)\end{aligned} \quad (1)$$

where $u_i(t)$ is the control input which is the desired acceleration for vehicle $i$ at time $t$; $\tau_i$ represents the inertial lag of the longitudinal dynamics of vehicle $i$. We consider the heterogeneity of vehicles where different vehicles have different values of the inertial lag $\tau_i$.

We assume each vehicle in the platoon can measure the distance $d_i$ and the front-vehicle velocity $v_{i-1}$. Combined with vehicle dynamics, a distributed cascade PID control algorithm is proposed to perform longitudinal control of the vehicle platoon. It is worth noting that the cascade PID is a double control structure with the inner and the outer loops. The outer loop PID controls the distance between adjacent vehicles, and the inner loop PID controls the vehicle speed. The output of the outer loop is used as the input of the inner loop and compared with the feedback value of the inner loop. This forms the whole inner-outer loop double-layer control. The control structure of each subsystem is shown in Fig. 3.

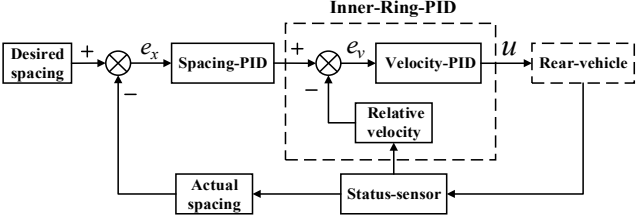

Fig.3. The cascade PID control structure

Let $d_i$ denote the measured distance between the front bumper of the following vehicle $i$ and the rear bumper of the leading vehicle $i-1$, we can obtain:
$$d_i(k) = x_{i-1}(k) - x_i(k) - l_{i-1} \quad (2)$$
where $k$ denotes the number of sampling period; $l_{i-1}$ denotes the length of the $(i-1)th$ vehicle. The desired spacing control strategy we adopt here is the constant time headway spacing strategy. The desired relative distance $S_i$ from the predecessor $i$-1's rear bumper to the follower $i$'s front bumper is defined as:
$$S_i(k) = d_0 + v_i(k)h_t \quad (3)$$
where $d_0$ and $h_t$ represent the minimum safe distance and the constant headway time, respectively. To simplify the problem, we assume all subsystems have the same $d_0$ and $h_t$.

The goal of the cascade PID control algorithm is to help each subsystem maintain the desired spacing and the consistent speed, and respond to any vehicle disturbance quickly such that the stability of the whole platoon can be maintained. The spacing error $e_{xi}$ and the velocity error $\dot{e}_{xi}$ which are used to measure the control target of the subsystem $P_{i,i-1}$ are defined as:
$$e_{xi}(k) = d_i(k) - S_i(k) \quad (4)$$
$$\dot{e}_{xi}(k) = v_{i-1}(k) - v_i(k) \quad (5)$$

The outer loop PID control equation of the subsystem $P_{i,i-1}$ is given by:
$$\partial_i(k) = K_{px}^i e_{xi}(k) + K_{ix}^i \sum_{j=0}^{k} e_{xi}(j) + K_{dx}^i [e_{xi}(k) - e_{xi}(k-1)] \quad (6)$$
where $\partial_i$ denotes the output of the outer loop PID; $K_{px}^i$, $K_{ix}^i$, $K_{dx}^i$ are the parameters of the outer loop PID.

The inner loop PID control equation of the subsystem $P_{i,i-1}$ is given by:
$$u_i(k) = K_{pv}^i e_{vi}(k) + K_{iv}^i \sum_{j=0}^{k} e_{vi}(j) + K_{dv}^i [e_{vi}(k) - e_{vi}(k-1)] \quad (7)$$
$$e_{vi}(k) = \partial_i(k) - \dot{e}_{xi}(k) \quad (8)$$
where $e_{vi}$ and $u_i$ represent the input and the output of the inner loop PID, respectively; $K_{pv}^i$, $K_{iv}^i$, $K_{dv}^i$ are the parameters of the inner loop PID.

Equation (1) can be discretized by using the method of difference approximation, and the discrete model of the acceleration is obtained as:
$$a_i(k+1) = \left(1 - \frac{T_s}{\tau_i}\right) a_i(k) + \frac{T_s}{\tau_i} u_i(k) \quad (9)$$
with constraints:
$$\begin{aligned} u_{min} \leq u_i(k) \leq u_{max} \\ a_{min} \leq a_i(k) \leq a_{max} \\ v_{min} \leq v_i(k) \leq v_{max} \end{aligned} \quad (10)$$
where $T_s$ represents the sampling time; $u_{min}$ and $u_{max}$ represent the minimum and the maximum desired acceleration; $a_{min}$ and $a_{max}$ are the minimum and the maximum acceleration; $v_{min}$ and $v_{max}$ are the minimum and the maximum velocity allowed on the road. The first constraint considers the vehicle performance limitation; the second constraint relates to passenger comfort and the third constraint reflects the road condition.

*C. Cooperative lane-changing model*

Fig. 4 shows the flow chart of the proposed cooperative lane-changing model which includes three parts: 1) lane-changing decision; 2) dynamic trajectory planning; 3) trajectory tracking.

*1) Lane-changing decision*

The status of the platoon vehicles can be divided into three types: 1) following the platoon, 2) joining the platoon and 3) leaving the platoon. This section discusses the second type of behavior. If a vehicle wants to merge into the platoon on the target lane, the vehicle and the platoon make cooperative decisions to complete the maneuver as shown in Fig. 5. $FV1$ and $FV2$ are the following vehicles in the original platoon on the target lane.

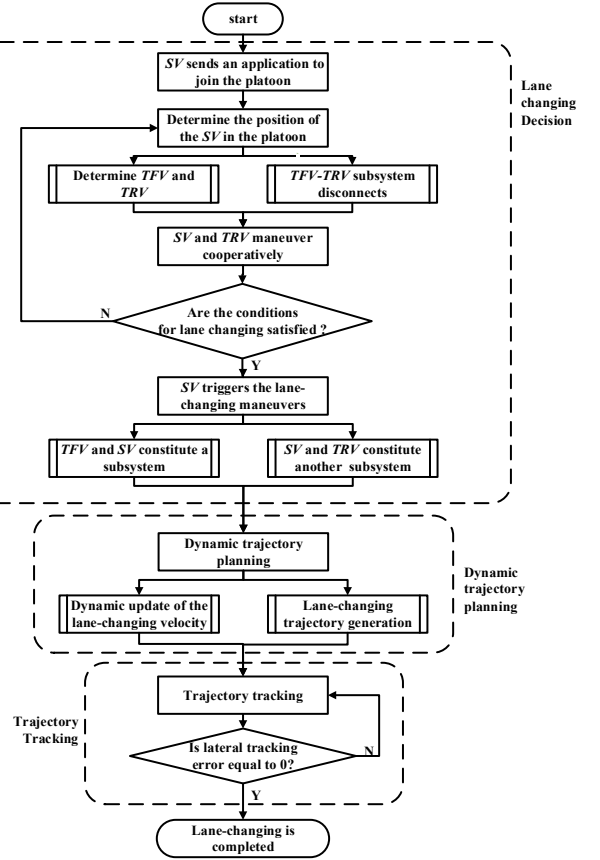

Fig. 4. The flow chart for dynamic lane changing model

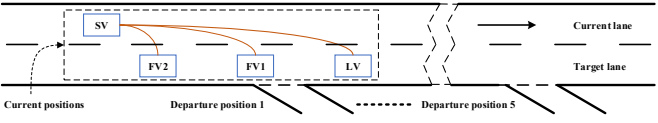

Fig. 5. Cooperative lane-change illustration

Firstly, the position of the SV in the platoon is determined according to the following assumption and rules defined as:
1. The location information of each vehicle where it leaves the freeway is available beforehand.
2. The distance between the departure position and the current position is used to determine the relative position of each vehicle in the platoon: the vehicle with the smallest distance is at the end of the platoon, while the vehicle with the largest distance is at the head of the platoon. That is to say, the position sequence of each vehicle in the platoon is determined according to the order in which it leaves the platoon. The last vehicle in the platoon always leaves the platoon first which reduces the impact on the stability of the whole platoon system.

Meanwhile, the *FV1* and the *FV2* become the *TFV* and the *TRV*, respectively, as shown in Fig. 6. The *SV* sends the request of joining the platoon to the *TFV* and the *TRV*.

Secondly, the subsystem composed of the *TRV* and the *TFV* is disconnected. The *TFV* and the vehicles in front of the *TFV* continue to move with the current status. We define the position of the *SV* at the moment when Equation (11) is satisfied as the starting point of lane-changing. The *TRV* and the *SV* maneuver cooperatively to meet the conditions of the lane-changing:1) The *SV* accelerates appropriately when the *SV* is behind the starting point of lane-changing, or decelerates appropriately when the *SV* is in front of the starting point of lane-changing; 2) The *TRV* decelerates to provide a safe lane-changing space for the *SV*. Let $d_{SV}$ and $d_{TRV}$ denote the longitudinal spacing between the *SV* and the *TFV*, and the longitudinal spacing between the *TRV* and the *SV*, respectively. Lane-changing is performed when the two spacings ($d_{SV}$, $d_{TRV}$) satisfy the following conditions:

$$\begin{matrix} d_{SV} = S_{SV} \\ d_{TRV} \geq d_0 \end{matrix} \quad (11)$$

where $S_{SV}$ represents the desired relative distance of the *SV*; $d_0$ represents the minimum safety distance.

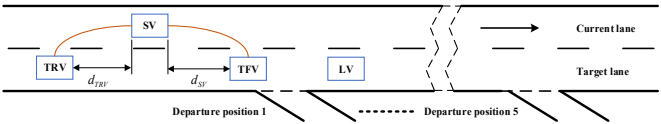

Fig.6. Lane-change spacing for vehicles

Thirdly, when the lane-changing is triggered, the *SV* and the *TFV* constitute a new subsystem. Meanwhile, the *TRV* and the *SV* constitute another subsystem to ensure safe longitudinal driving. Then, the lane-changing maneuver is carried out according to the model proposed in subsection 2) in section II until the whole lane-changing process is completed.

*2) Dynamic trajectory planning*

In this subsection, we integrate the distributed cascade PID control algorithm with the sine function [37] into a novel dynamic lane-change trajectory planning model. The DCPID model described by Equations (1-9) is applied to update the acceleration of the *SV* during the lane-changing process, where we define a function $F_{DCPID}$ as shown in Equation (12). The sine function is used to derive the lane-changing trajectory. Let $x_{SV}^0$, $y_{SV}^0$ and $y_{TFV}^0$ denote the lateral and the longitudinal position of vehicle *SV* and the position of the *TFV* at the beginning of the lane-changing process, respectively. The dynamic lane-changing trajectory planning part including dynamic update of the lane-changing speed and the lane-changing trajectory generation is given by:

$$\begin{cases} y_d^0 = y_{TFV}^0 - y_{SV}^0 \\ a_{SV}(t+T_s) = F_{DCPID}\big(x_{SV}(t), x_{TFV}(t), v_{SV}(t), v_{TFV}(t)\big) \\ v_{SV}(t+T_s) = v_{SV}(t) + a_{SV}(t+T_s) * T_s, \ t \in [t_0, t_e] \\ y_r(x_r) = y_{SV}^0 + \frac{y_d^0}{2\pi}\left\{\frac{2\pi}{M}(x_r - x_{SV}^0) - \sin\left[\frac{2\pi}{M}(x_r - x_{SV}^0)\right]\right\} \\ M = v_{SV}(t)\sqrt{\frac{2|y_d^0|}{a_p}} \\ x_r \in [x_{SV}^0, x_{SV}^0 + M] \end{cases} (12)$$

where $y_d^0$ is the lateral distance at the beginning of lane-changing maneuver between the *TFV* and the *SV*; $(x_r, y_r)$ is the desired location of vehicle *SV*; $t_0$ and $t_e$ represent the start time and the end time of lane-changing, respectively; $v_{SV}$ represents the velocity of the *SV*, which is the safe lane-changing speed derived based on the DCPID longitudinal control algorithm; $a_p$ is the planned acceleration considering the comfort of lane-changing. It is worth noting that $a_p$ is a parameter which has physical significance. Thus, we can obtain the appropriate value of $a_p$ by sensitivity analysis. $v_{SV}$ and $a_p$ jointly determine the safety, real-time capability and ride comfort of the whole lane-changing process.

Let $y_r{}'$, $y_r{}''$ and $K$ denote the first derivative of $y_r$, the second derivative of $y_r$ and the curvature of $y_r$ respectively, as given by:

$$y_r{}'(x_r) = \frac{y_d^0}{M}\left\{1 - \cos\left[\frac{2\pi}{M}(x_r - x_{SV}^0)\right]\right\} \quad (13)$$

$$y_r{}''(x_r) = \frac{\pi a_p}{v_{SV}{}^2(t)}\sin\left[\frac{2\pi}{M}(x_r - x_{SV}^0)\right] \quad (14)$$

$$K = \frac{y_r{}''}{(1+y_r{}'^2)^{3/2}} \quad (15)$$

Let $\varphi_r$ and $\delta_{fr}$ denote the desired yaw angle and the desired front wheel angle of *SV*, which can be calculated as:

$$\varphi_r = \tan^{-1}(y_r{}') \quad (16)$$
$$\delta_{fr} = \tan^{-1}(L * K) \quad (17)$$

where $L$ denotes the distance between the front and the rear axles; $\varphi_r$ and $\delta_{fr}$ are in radians.

After the lane-changing decision, vehicle *SV* and vehicle *TFV* are regarded as a new subsystem controlled by the DCPID control algorithm. Under the premise of avoiding collision between the two vehicles, the DCPID control algorithm can calculate the acceleration of the lane-changing vehicle and update the velocity of the vehicle in real time by considering the relative speed and distance between the two vehicles. Then, the speed is transferred to the sine function model instantly to plan a lane-changing path. The proposed model can update the lane-changing path in real time according to the speed change of the *TFV*, which not only ensures the safety of lane-changing, but also improves the flexibility of lane-changing.

*3) Trajectory tracking*

In order to simplify the microscopic vehicle control, we use a kinematic model of three degrees of freedom to describe the vehicle state. The formula is given by:

$$[\dot{x} \quad \dot{y} \quad \dot{\varphi}]^T = [\cos\varphi \quad \sin\varphi \quad \tan\delta_f/l]^T \cdot v \quad (18)$$

where $(x, y)$ is the coordinates of the rear axle center of the vehicle; $\varphi, \delta_f, l, v$ are the yaw angle, the front wheel angle, distance between the front and the rear axles, and the velocity, respectively.

For the reference trajectory generated in subsection II.C. 2), we adopt the model predictive control (MPC) algorithm to perform the real-time tracking control as given by:

$$\widetilde{\Phi}(k+1) = A(k)\widetilde{\Phi}(k) + B(k)\Delta\tilde{u}(k) \quad (19)$$

with

$$\widetilde{\Phi}(k) \triangleq \begin{bmatrix} X(k) - X_r(k) \\ u(k) - u_r(k) \end{bmatrix} \quad (20)$$

where $X = [x \quad y \quad \varphi]^T$ and $u = [v \quad \delta_f]^T$ are the current state and the current control variables, respectively; $X_r = [x_r \quad y_r \quad \varphi_r]^T$ and $u_r = [v_r \quad \delta_{fr}]^T$ are the desired state and the desired control variables respectively obtained from the reference trajectory.

We design the objective function to ensure that the vehicle can track the reference trajectory quickly and smoothly. The minimum cost function of the MPC is given by:

$$\min J(k) = [\Delta U(k)^T, \varepsilon]^T H(k)[\Delta U(k)^T, \varepsilon] \\ + f^T(k)[\Delta U(k)^T, \varepsilon] \quad (21)$$

$$H(k) \triangleq \begin{bmatrix} \bar{B}_k^T \bar{Q}\bar{B}_k + \bar{R} & 0 \\ 0 & \rho \end{bmatrix} \quad (22)$$

$$f(k)^T = [2\bar{A}_k^T \bar{Q}\bar{B}_k \quad 0]$$

with constraints:

$$\begin{cases} \Delta U_{min} \leq \Delta U(k) \leq \Delta U_{max} \\ U_{min} \leq A\Delta U(k) + U(k-1) \leq U_{max} \end{cases} \quad (23)$$

where $k$ is the current time step; $\Delta U(k)$ is the system control variable increment at time step $k$; $U_{min}, U_{max}$ are the limits of the control variable; $\Delta U_{min}, \Delta U_{max}$ are the limits of the control variable increment; $\varepsilon$ is the relaxation factor; $\rho$ is the weight coefficient; $\bar{Q}, \bar{R}$ are the weight matrices; $\bar{A}_k, \bar{B}_k$ are the prediction matrices; $A$ is the coefficient matrix of the constraint equation; $H(k)$ and $f(k)$ are the coefficient matrices of the standard form of the quadratic programming problem.

Therefore, this problem is transformed into a standard Quadratic Programming (QP) problem under the MPC framework. The first part of Equation (21) reflects the tracking ability of the control system to the reference trajectory, while the second part represents the constraints of the control variables of the system.

## III. STABILITY ANALYSIS OF THE DCPID ALGORITHM

This section presents the stability analysis of the proposed DCPID algorithm, which includes local stability analysis and asymptotic stability analysis. The vehicle can maintain an equilibrium state (for example, the desired speed and the desired spacing) under disturbance, which is called local stability. The amplitude of the disturbance (for example, deviation from the desired speed and the desired spacing) gradually attenuates as it propagates from downstream to upstream in the vehicle platoon, which is called asymptotical stability [17]. We use the stability analysis method proposed by [38] to derive sufficient conditions for local stability and asymptotical stability. These conditions show that the stability can be achieved by proper parameter tuning of the proposed DCPID algorithm.

### A. Derivation of Partial Differential Equation

The general formation of the DCPID longitudinal control can be expressed as:

$$\dot{v}_i(t) = f(v_i(t), \dot{e}_{xi}(t), d_i(t)) \quad (24)$$

where $\dot{e}_{xi}(t) = v_{i-1}(t) - v_i(t)$ represents the speed difference between the preceding vehicle *i-1* and the following vehicle *i* at time *t*; $d_i(t)$ represents the distance between the following vehicle *i* and the preceding vehicle *i-1* at time *t*; $\dot{v}_i(t)$ is the acceleration of vehicle *i* at time *t*.

The partial differential values $f_v$, $f_{\dot{e}_x}$ and $f_d$ of equation (24) for $v_i(t)$, $\dot{e}_{xi}(t)$ and $d_i(t)$ in the equilibrium state $(v_e, 0, d_e)$ can be obtained as:

$$\begin{cases} f_v = \frac{\partial f(v_i(t),\dot{e}_{xi}(t),d_i(t))}{\partial v_i(t)}|_{(v_e,0,d_e)} \\ f_{\dot{e}_x} = \frac{\partial f(v_i(t),\dot{e}_{xi}(t),d_i(t))}{\partial \dot{e}_{xi}(t)}|_{(v_e,0,d_e)} \\ f_d = \frac{\partial f(v_i(t),\dot{e}_{xi}(t),d_i(t))}{\partial d_i(t)}|_{(v_e,0,d_e)} \end{cases} \quad (25)$$

By substituting equations (1)-(9) into equation (25), we can obtain:

$$\begin{cases} f_v = -\frac{T_s h_t}{\tau_i} \left[\frac{1}{2}K_{ix}^i K_{iv}^i t^2 + (K_{ix}^i K_{pv}^i + K_{px}^i K_{iv}^i)t + K_{px}^i K_{pv}^i + K_{ix}^i K_{dv}^i\right] \\ f_{\dot{e}_x} = -\frac{T_s}{\tau_i}(K_{iv}^i t + K_{pv}^i) \\ f_d = \frac{T_s}{\tau_i}(K_{ix}^i t + K_{px}^i) \end{cases} \quad (26)$$

### B. Sufficient Condition of DCPID stability

The conditions of local stability and asymptotic stability of the DCPID are stated as follows.

*Theorem 1:* The vehicle controlled by the DCPID is locally stable if satisfying:

$$f_v - f_{\dot{e}_x} < 0 \quad (27)$$

*Theorem 2:* The vehicle platoon under the DCPID is asymptotically stable if satisfying:

$$\frac{1}{2}(f_v)^2 - f_v f_{\dot{e}_x} - f_d > 0 \quad (28)$$

Equations (27) and (28) show that the stability of the DCPID is closely related to its parameters. Therefore, we can ensure the stability of the DCPID through proper parameter tuning. This provides theoretical support for the parameter settings in the simulation experiment.

## IV. SENSITIVITY ANALYSIS OF PARAMETERS

### A. Sensitivity analysis for DCPID control parameters

In the DCPID model, six key control parameters need to be determined including 3 parameters ($K_{pv}^i$, $K_{iv}^i$, $K_{dv}^i$) for the inner loop control and 3 parameters ($K_{px}^i$, $K_{ix}^i$, $K_{dx}^i$) for the outer loop control as indicated in Equations (6) and (7). We denote $P =$

$[K_{px}\ K_{ix}\ K_{dx}\ K_{pv}\ K_{iv}\ K_{dv}]$ as the parameter set. Different parameter sets have significant influence on the performance of platoon control. Thus, it is necessary to explore the appropriate parameter range for different platoon states. A two-vehicle system is used to analyze the sensitivity of DCPID control parameters. We define the system state $E_s$ as:

$$E_s = [ex\ \ ev] \quad (29)$$

with:

$$\begin{aligned} ex &= d - S \\ ev &= v_1 - v_2 \end{aligned} \quad (30)$$

where $d$ and $S$ represent the actual distance and the desired distance between the two vehicles, respectively; $v_1$ and $v_2$ represent the velocities of the front vehicle and the rear vehicle, respectively; $ex$ and $ev$ represent the spacing error and the velocity error of the two-vehicle system. In addition, the negative values of $ex$ and $ev$ indicate the rear vehicle needs to slow down, and the positive values of $ex$ and $ev$ indicate that the rear vehicle needs to speed up.

In order to evaluate the effectiveness of the proposed control algorithm, two performance indexes, namely, the steady state adjustment time $t_{steady}$ and the velocity overshoot rate $\eta$ are chosen. The latter is calculated as:

$$\eta = \frac{v_{overshoot}}{v_1} \times 100\% \quad (31)$$

where $v_{overshoot}$ is the velocity overshoot; $v_1$ is the steady state velocity of the system.

In the following, we explore the relationship between the DCPID parameter values and $E_s$ in different states. According to the values of $E_s$, we can roughly divide the states of each subsystem into four cases:

**Case 1**: the subtle disturbance caused by the deceleration of the leading vehicle under the steady state;

**Case 2**: The unstable state of the system caused by the existence of $ev$;

**Case 3**: The unstable state of the system caused by the existence of $ex$;

**Case 4**: The unstable state of the system caused by the existence of both $ev$ and $ex$.

Fig. 7 shows the performance of the system in terms of the velocity and the spacing errors under different parameter sets P for case 4 with $E_s = [5\ m, -2\ m/s]$. In Fig. 7 (a), the velocity overshoot phenomenon can be observed for the six parameter sets indicated with solid lines. It can also be seen from Fig. 7 (b) that the system takes longer steady-state adjustment time or even cannot reach the steady state with those parameter sets (solid lines). The system performs much better in terms of $t_{steady}$ and $v_{overshoot}$ with parameter sets indicated with three dashed lines, among which the parameter set $P = [5\ 0\ 10\ 5\ 0\ 0]$ outperforms.

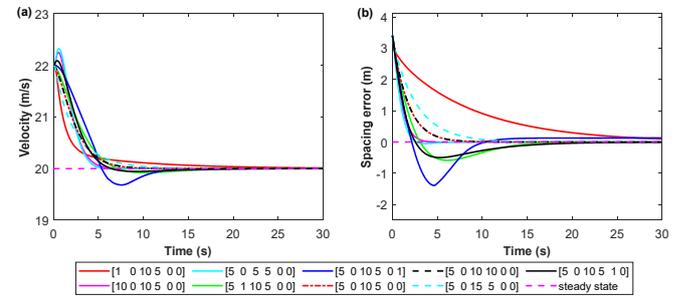

Fig.7. Comparison of the performance of different parameter sets for case 4 with $E_s = [5\ m, -2\ m/s]$: (a) velocity; (b) spacing error

TABLE I
PARAMETER ANALYSIS FOR DCPID

| Cases | Scenarios | $ev$ (m/s) | $ex$ (m) | $K_{px}$ | $K_{ix}$ | $K_{dx}$ | $K_{pv}$ | $K_{iv}$ | $K_{dv}$ | $\eta$ | $t_{steady}$ |
|---|---|---|---|---|---|---|---|---|---|---|---|
| Case 1: small disturbance | 1 | 0 | 0 | 300 | 0 | 0 | 8 | 0 | 2 | 0 | 6 |
| Case 2: $ev$ is larger | 2 | 1 | 0 | 10 | 0 | 1-5 | 1-30 | 0 | 0 | 0 | 4 |
| | 3 | 2 | 0 | 10 | 0 | 1-3 | 1-30 | 0 | 0 | 0 | 4 |
| | 4 | 3 | 0 | 10 | 0 | 1-3 | 1-30 | 0 | 0 | 1.5 | 5 |
| | 5 | 5 | 0 | 5-15 | 0 | 3-5 | 1-10 | 0 | 0 | 6.0 | 8 |
| | 6 | 7 | 0 | 5-15 | 0 | 5-10 | 1-10 | 0 | 0 | 8.0 | 10 |
| | 7 | 10 | 0 | 5-15 | 0 | 9-10 | 1-8 | 0 | 0 | 13.3 | 15 |
| Case 3: $ex$ is larger | 8 | 0 | -10 | 8-15 | 0 | 8-10 | 1-8 | 0 | 0 | 0 | 10 |
| | 9 | 0 | -5 | 8 | 0 | 10 | 5 | 0 | 0 | 0 | 15 |
| | 10 | 0 | -2 | 8 | 0 | 10 | 5 | 0 | 0 | 0 | 15 |
| | 11 | 0 | 1 | 8 | 0 | 10 | 5 | 0 | 0 | 0 | 15 |
| | 12 | 0 | 3 | 8 | 0 | 10 | 5 | 0 | 0 | 0 | 15 |
| | 13 | 0 | 7 | 8 | 0 | 10 | 5 | 0 | 0 | 0 | 15 |
| | 14 | 0 | 10 | 8-15 | 0 | 8-10 | 1-8 | 0 | 0 | 0 | 16 |
| Case 4: $ev$ and $ex$ are both large | 15 | 5 | 5 | 10 | 0 | 15 | 2 | 0 | 0 | 12.0 | 15 |
| | 16 | -5 | -5 | 2 | 0 | 5 | 2 | 0 | 0 | 4.0 | 10 |
| | 17 | -5 | 5 | 10 | 0 | 5-15 | 1-8 | 0 | 0 | 15.0 | 15 |
| | 18 | 5 | 2 | 5 | 0 | 10 | 5 | 0 | 0 | 4.0 | 15 |
| | 19 | 5 | -2 | 5 | 0 | 5 | 5 | 0 | 0 | 1.2 | 10 |
| | 20 | -5 | 2 | 5 | 0 | 10 | 5 | 0 | 0 | 6.0 | 15 |
| | 21 | -2 | 5 | 5 | 0 | 10 | 5 | 0 | 0 | 0 | 15 |
| | 22 | 2 | 5 | 5 | 0 | 10 | 5 | 0 | 0 | 6.8 | 20 |

We conduct parameter tuning for each case with different scenarios considering $t_{steady}$ and $v_{overshoot}$. The results of the control parameter sets are shown in Table I. As can be seen from the table, when only the spacing error $ex$ exists, no velocity overshoot occurs in the system with $\eta = 0$. However, when the velocity error $ev$ exists and is greater than $3m/s$, the velocity overshoot phenomenon occurs. This indicates that the DCPID algorithm gives better response to spacing than velocity. In addition, it can be clearly seen from Table I that the control parameters for case 1 are quite different from those for cases 2-4 which have similar parameters (or parameter range). Therefore, the values of control parameters are divided into two groups:

Group 1: control parameter set for case 1
Group 2: control parameter set for cases 2-4

In cases 2-4, the parameter $K_{px}$ of each scenario is around 8. The values of $K_{dx}$ and $K_{pv}$ are around 10 and 5, respectively. Among all cases, the values of $K_{ix}$ and $K_{iv}$ are 0. Thus, the final control parameter values for these two groups are shown in Table II.

TABLE II
FINAL CONTROL PARAMETER VALUES FOR DCPID

| Group | $K_{px}$ | $K_{ix}$ | $K_{dx}$ | $K_{pv}$ | $K_{iv}$ | $K_{dv}$ |
|---|---|---|---|---|---|---|
| 1 | 300 | 0 | 0 | 8 | 0 | 2 |
| 2 | 8 | 0 | 10 | 5 | 0 | 0 |

## B. Feasible domain for the parameter $a_p$ of trajectory planning

In the cooperative lane-changing model, the feasible domain of $a_p$ in Equation (12) is determined by considering vehicle ride comfort. We define the yaw-rate $\omega$ of vehicle as an indicator to measure the vehicle ride comfort.

$$\omega(t) = \frac{a_y(t)}{v_{SV}(t)} \quad (32)$$

where $\omega(t)$, $a_y(t)$ and $v_{SV}(t)$ represent the yaw-rate, the lateral acceleration and the velocity of the $SV$ at time $t$, respectively. The upper bound of the yaw-rate $\omega_{upper\_bound}$ that can guarantee comfort during lane-changing process can be described by [33]:

$$\omega_{upper\_bound} = 0.85 \frac{a_{ymax}}{v_{SV}} = \frac{0.425}{v_{SV}} (rad/s) \quad (33)$$

where $a_{ymax}$ represents the maximum lateral acceleration and $a_{ymax} = 0.5\ m/s^2$ [33], [39]. Ride comfort can be guaranteed during the lane-changing process when the following conditions are met:

$$\begin{aligned} |\omega|_{max} &\leq \omega_{upper\_bound} \\ |\omega|_{max} &= \max_{t_0 \leq t \leq t_e} |\omega(t)| \end{aligned} \quad (34)$$

Fig. 8 shows the impact of $a_p$ on the yaw-rate at three different speeds. The three dashed lines correspond to $\omega_{upper\_bound}$ of the $SV$ at $20m/s$, $25m/s$ and $30m/s$ from top to bottom, respectively. The part below the dashed line indicates that the corresponding $a_p$ value can ensure the requirement of riding comfort, while the part above the dashed line indicates that the corresponding $a_p$ value cannot meet the requirement. It can also be seen from Fig. 8 that the feasible domain of $a_p$ satisfying ride comfort decreases as $v_{SV}$ increases. When $a_p$ is determined, the difference between $|\omega|_{max}$ and $\omega_{upper\_bound}$ decreases as $v_{SV}$ increases. This indicates that the requirements for ride comfort become higher as the velocity increases during the lane-changing process. The $\omega_{upper\_bound}$ values and the feasible domains of $a_p$ under three velocities of the $SV$ are shown in Table III.

TABLE III
THE VALUES OF $\omega_{upper\_bound}$ AND THE FEASIBLE DOMAINS OF $a_p$ AT DIFFERENT VELOCITIES OF THE $SV$

| $v_{SV}(m/s)$ | $\omega_{upper_{bound}}(rad/s)$ | Feasible domain of $a_p$ |
|---|---|---|
| 20 | 0.0212 | (0, 0.122] |
| 25 | 0.0170 | (0, 0.114] |
| 30 | 0.0142 | (0, 0.106] |

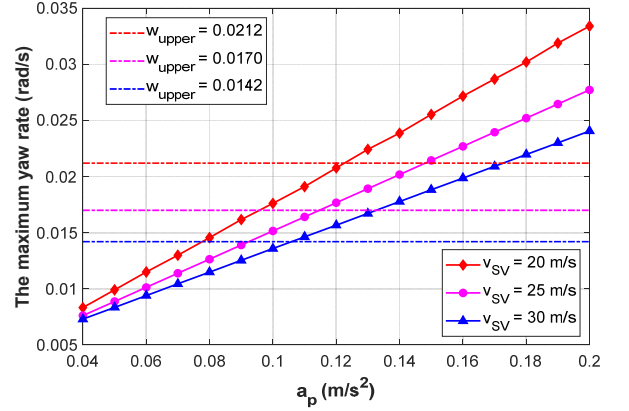

Fig.8. Effect of $a_p$ on $|\omega|_{max}$ at three different speeds

## V. SIMULATIONS AND RESULTS ANALYSIS

### A. Longitudinal platooning without lane-changing

To demonstrate the performance of the proposed DCPID model, a numerical simulation is conducted considering a heterogeneous platoon composed of one leading vehicle and seven following vehicles.

The setting of the main simulation parameters is shown in Table IV. Since the state of each subsystem is similar, we set the same DCPID parameters for each subsystem. The final parameters of the DCPID model are determined through the sensitivity analysis as discussed in Section IV with $K_{px}^i = 8, K_{ix}^i = 0, K_{dx}^i = 10, K_{pv}^i = 5, K_{iv}^i = 0, K_{dv}^i = 0$. In addition, the above parameters satisfy the conditions of local stability and asymptotic stability conditions in (27) and (28).

TABLE IV
THE MAIN SIMULATION PARAMETERS[19], [40]

| | Parameter | Notation | Unit | Value |
|---|---|---|---|---|
| | sampling time | $T_s$ | s | 0.02 |
| | time headway | $h_t$ | s | 0.8 |
| | minimum safe distance | $d_0$ | m | 4 |
| Control parameters | minimum control value | $u_{min}$ | $m/s^2$ | -3 |
| | maximum control value | $u_{max}$ | $m/s^2$ | 3 |
| | minimum acceleration | $a_{min}$ | $m/s^2$ | -3 |
| | maximum acceleration | $a_{max}$ | $m/s^2$ | 3 |
| | length of the vehicles | $l_i$ | m | 5 |
| Vehicle parameters | inertial lag of longitudinal dynamics | $\tau_i$, $i \in (2,8)$ | s | [0.51,0.75, 0.78,0.70, 0.73,0.72, 0.62] |

### 1) Numerical results

To demonstrate that the parameters of the DCPID obtained through stability analysis proposed in Section III and parameter sensitivity analysis proposed in Section IV can provide robust control for tracking and adjusting the desired spacing and speed both in the stable state and the unstable state, we use $ex$ and $ev$ defined in Equations (29) and (30) to describe the state of the vehicle platoon. We increase $ex$ linearly from $-10m$ to $10m$ with a step length of 1m and increase $ev$ linearly from $-5m/s$ to $5m/s$ with a step of $0.5m/s$. In total 400 scenarios are obtained. By simulating each scenario and calculating the steady-state adjustment time $t_{steady}$ and the velocity overshoot rate $\eta$ defined in Equation (31), we can evaluate whether the DCPID has good tracking and adjustment performance.

shows the trajectory in the longitudinal direction where no collision occurs during the entire process, which indicates that the safety of the DCPID controlled platoon can be assured.

We further analyze the simualtion results for all 400 scenarios and conclude that the DCPID can guarantee the string stability of the platoon for all scenarios. The overshoot rate $\eta$ of each scenario is shown in Fig. 10. From the figure, we can observe that the velocity overshoot rates of most scenarios are rather small with $\eta$ <5% (the dark blue area).

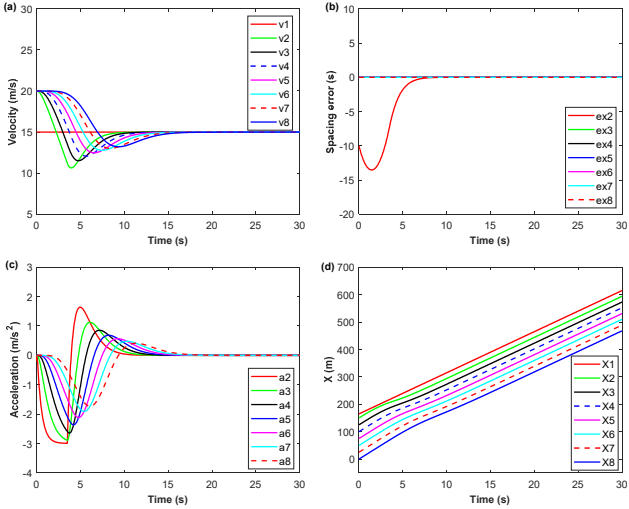

Fig. 9. Simulation results of the DCPID Platoon under extreme conditions with $ex = -10m$ and $ev = -5m/s$: (a) velocity; (b) spacing error; (c) acceleration; (d) the longitudinal trajectory

Fig. 9 shows the response of the following vehicles to the speed and the desired spacing under the extreme conditions where $ex = -10m$ and $ev = -5m/s$. We can observe that all vehicles can adjust the velocity and spacing to the stable state despite that the initial state of the platoon is relatively extreme. This also reveals the robustness of the DCPID model. Fig. 9 (d)

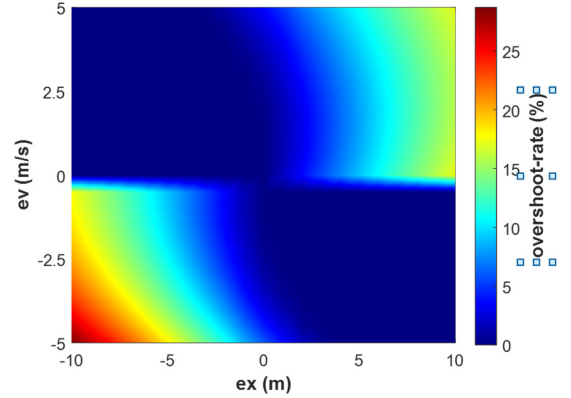

Fig. 10. The velocity overshoot rate with $ex \in [-10m, 10m]$ and $ev \in [-5m/s, 5m/s]$

*2) Comparison Analysis*

In this section, we discuss the performance comparison of our proposed DCPID with the DMPC proposed by [17] and the single PID. We consider a platoon of 8 vehicles using the same scenario setting as [17]. The initial spacing errors are set to be 2, 0.1, 0.1, 0.1, 0.1, 0.1, 0.1, and both the initial speed errors and acceleration errors are set to be 0. As can be seen from Fig. 11, all algorithms can effectively adjust the spacing error to zero. Nevertheless, the proposed DCPID provides better tracking performance with the smallest spacing error variation. Compared with a single-loop PID, the parallel adjustment of the inner and outer loops of the DCPID achieves a faster response speed and a smaller overshoot.

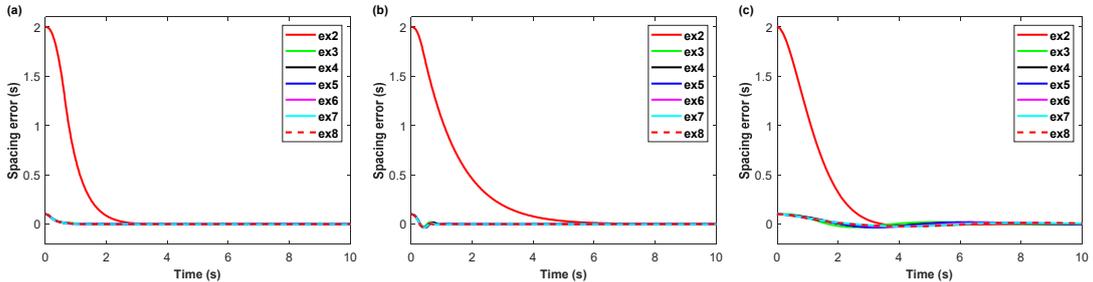

Fig. 11. Variation of spacing error under three methods: (a) DCPID; (b) single-PID; (c) DMPC in [17]

In order to further investigate anti-interference ability of the DCPID, we add a disturbance $\epsilon_u$, which is considered as the system internal interference, to the control variable $u$. We consider a platoon of 8 vehicles in a steady state (the initial spacing errors and speed errors are all 0). During the time between $t=6s$ and $t=8s$, interference $\epsilon_u = 3$ is added to $u$ of the leading vehicle. As can be seen from Fig. 12, all these three algorithms can well suppress the interference and restore the vehicle platoon to a stable state. We can also observe that with the proposed DCPID approach, the amplitudes of spacing error and speed error caused by interference are significantly smaller than those with the DMPC and the single PID. This indicates that the DCPID has stronger anti-interference ability. The underlying reason is that the double-layer structure of the

DCPID can adjust the spacing and the speed parallelly to smooth the disturbance.

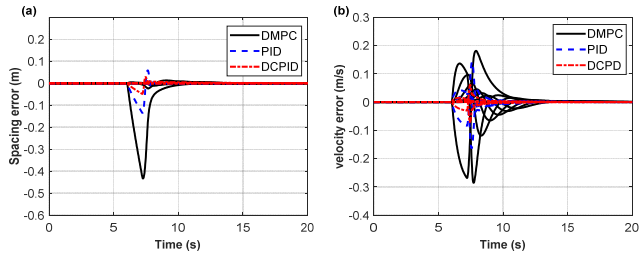

Fig. 12. Comparison analysis of the DCPID with respect to single-PID and DMPC in [17]: (a) spacing error; (b) velocity error

*B. Cooperative lane-changing*

To evaluate the performance of the dynamic cooperative lane-changing model, we conducted several simulation tests listed in Table V. The DCPID-controlled platoon is composed of four vehicles driving on the target lane and a SV driving on the adjacent lane. The initial state of each vehicle is defined as $\theta_i^0 = (x_{i0}, y_{i0}, v_{i0})$, $i \in \{1,2,3,4,SV\}$, where vehicles $i = \{2,3\}$ represent the *TFV* and the *TRV* in the platoon, respectively. We also consider different speed states of the *TFV* for different scenarios. The speed is set as constant for scenarios 1-3 and 5 and to be variable for scenario 4. The planned acceleration $a_p$ in the planning process is set to be $0.1 m/s^2$ according to Table III to ensure ride comfort in the following scenarios, and the distance between the front and the rear axles $L = 2.9m$. The control parameters and vehicle parameters are the same as listed in Table IV. Each subsystem can adjust the DCPID parameter values according to its own state referring to Table II to ensure the optimal control of the whole system.

TABLE V
PARAMETER SETTINGS FOR DIFFERENT SCENARIOS

|  | Scenario 1 | Scenario 2 | Scenario 3 | Scenario 4 | Scenario 5 |
|---|---|---|---|---|---|
| $\theta_1^0$ | $(75, -1.875, 20)$ | $(75, -1.875, 20)$ | $(116, -1.875, 25)$ | $(116, -1.875, 25)$ | $(132, -1.875, 30)$ |
| $\theta_2^0$ | $(50, -1.875, 20)$ | $(50, -1.875, 20)$ | $(87, -1.875, 25)$ | $(87, -1.875, 25)$ | $(99, -1.875, 30)$ |
| $\theta_3^0$ | $(25, -1.875, 20)$ | $(25, -1.875, 20)$ | $(35, -1.875, 25)$ | $(29, -1.875, 25)$ | $(33, -1.875, 30)$ |
| $\theta_4^0$ | $(0, -1.875, 20)$ | $(0, -1.875, 20)$ | $(6, -1.875, 25)$ | $(0, -1.875, 25)$ | $(0, -1.875, 30)$ |
| $\theta_{SV}^0$ | $(10, 1.875, 20)$ | $(35, 1.875, 20)$ | $(50, 1.875, 25)$ | $(58, 1.875, 25)$ | $(66, 1.875, 30)$ |
| *TFV* speed | Constant | Constant | Constant | Variable | Constant |

Fig. 13 illustrates the initial position of vehicles in different scenarios. Fig. 13 (a)(c) show scenarios 1 and 3 with $d_{SV} > S_{SV}$. According to Equation (11), we can tell that the *SV* needs to accelerate to satisfy lane changing requirements; Fig. 13 (b) shows scenario 2 with $d_{SV} < S_{SV}$, indicating that the *SV* needs to decelerate to satisfy lane-changing requirements; Fig. 13 (e)(d) show scenarios 4 and 5 with $d_{SV} = S_{SV}$. In this case, the *SV* only needs to consider the longitudinal spacing between the *TRV* and itself ($d_{TRV}$) ensuring $d_{TRV} \geq d_0$, such that the *SV* can trigger the lane changing maneuver.

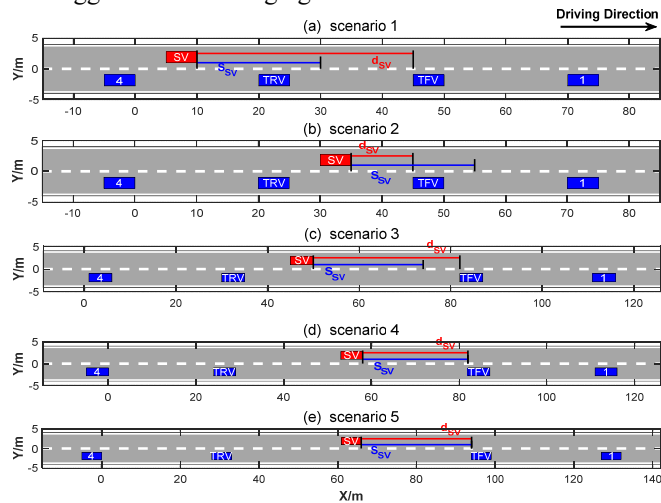

Fig. 13. The initial position of vehicles according to Table V' scenarios

Fig. 14(a) shows the speed dynamics for each vehicle during the cooperative lane-changing process in scenario 1. After accelerating and decelerating, the *SV* finally reaches the same speed as the preceding vehicles 1 and *TFV* in the platooning state. As for the *TRV* and vehicle 4, they need to decelerate at first to make space for the *SV* to change lanes and then accelerate to catch up the *SV*. The whole cooperative lane-changing process lasts about $16.36s$. After that, all vehicles drive with the same speed in a steady state. The lane changing process is illustrated by a three-dimensional plot of vehicle trajectories in Fig. 14(b). The velocity error and the spacing error dynamics are provided in Fig. 14 (c)(d), respectively.

Fig. 15 shows the performance of the *SV* during the lane-changing process in scenario 1. The comparison of the desired trajectory and the trajectory derived from the proposed model during the entire lane-change process is illustrated in Fig. 15 (a). It can be observed that the model derived trajectory is quite smooth and perfectly fit with the desired trajectory. Fig. 15 (b) shows the lateral tracking error with the maximum error of $0.001m$. This indicates that the proposed method can provide accurate trajectory tracking. The front wheel steering angle and the yaw angle of the *SV* are shown in Fig.15 (c) and Fig. 15 (d), respectively, which illustrate that the lane-changing is performed exactly as the desired.

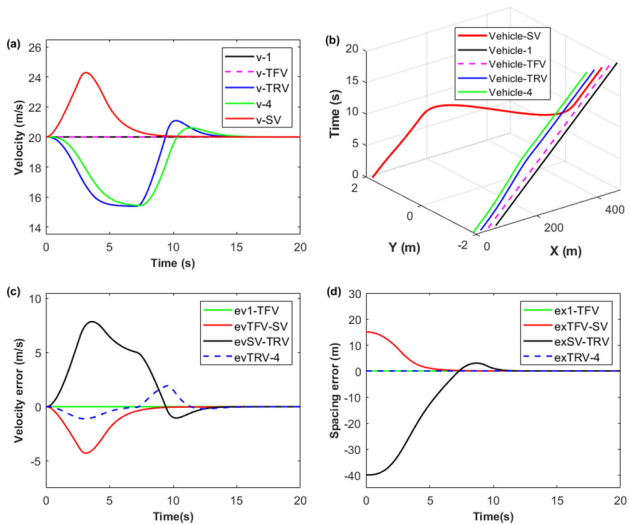

Fig. 14. Performance of the cooperative lane-changing in scenario 1: (a) velocity; (b) trajectory; (c) velocity error and (d) spacing error

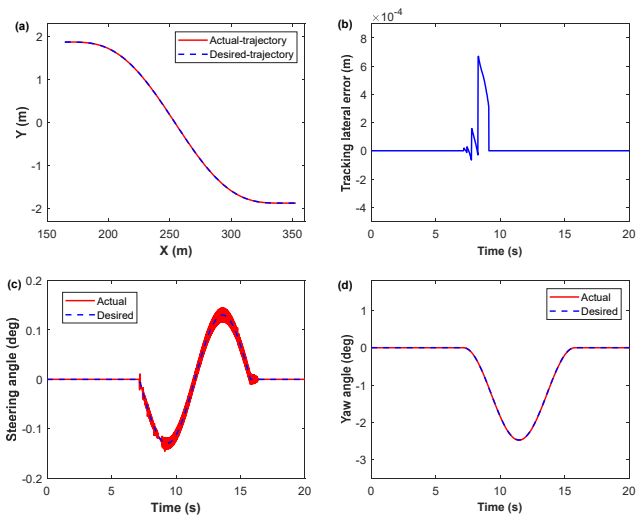

Fig. 15. Performance of the $SV$ during the lane-changing process in scenario 1: (a) trajectory; (b) lateral tracking error; (c) steering angle; (d) yaw angle

Fig. 16 shows the lateral tracking errors of the five scenarios. It can be seen that the lateral tracking errors of the $SV$ reach 0, indicating that the $SV$ can complete lane-change safely in all five scenarios. The maximum lateral tracking errors of scenarios 1, 2, and 3 are quite small about 0.001m; the maximum errors of scenarios 4 and 5 are slightly larger about 0.018m and 0.012m, respectively. The reasons behind are that the velocity of the $TFV$ is changing during the lane-changing process in scenario 4, causing the $SV$ to adjust its own velocity to avoid potential collisions. The fluctuation of the velocity also results in larger difference between the trajectories planned at different times, leading to the increased uncertainty of the lateral movement of the $SV$. In addition, the slight fluctuations may give rise to larger errors when the velocity is higher (e.g., scenario 5). This is also in line with the analysis in subsection IV.B where the lateral error increases and ride comfort decreases as the speed of lane-changing increases.

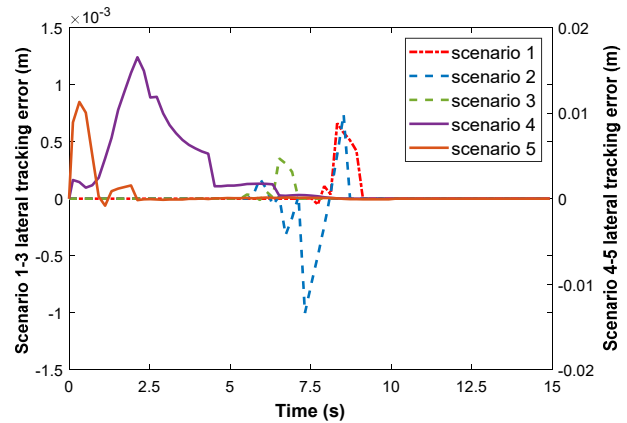

Fig. 16. The lateral tracking errors for five scenarios

Fig. 17 shows the longitudinal spacing errors of five scenarios. As can be seen the longitudinal spacing errors of all subsystems can return to the steady state without error in the five scenarios. This indicates that the $SV$ can successfully join the platoon and form a new steady platoon in all five scenarios.

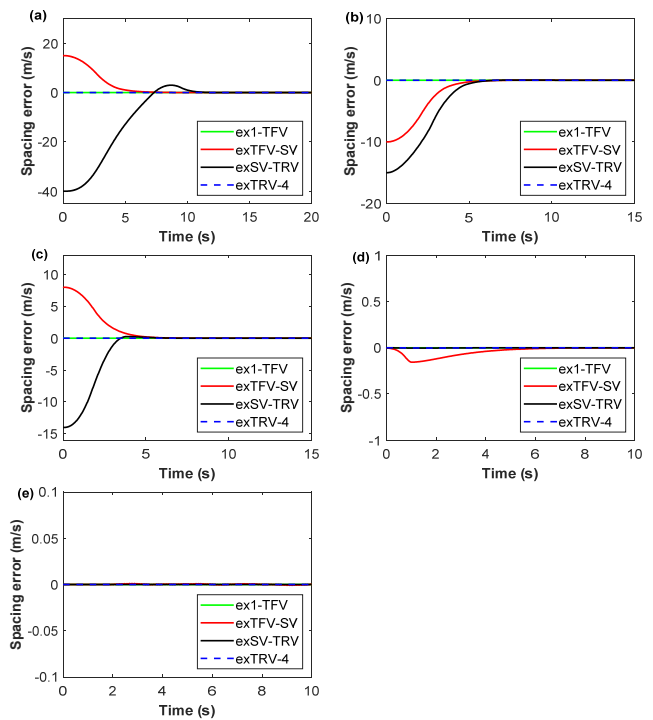

Fig. 17 The longitudinal spacing errors for five scenarios: (a) scenario 1; (b) scenario 2; (c) scenario 3; (d) scenario 4; (e) scenario 5

TABLE VI
THE START AND THE END TIME OF LANE-CHANGING IN FIVE SCENARIOS

| $t(s)$ | Scenario 1 | Scenario 2 | Scenario 3 | Scenario 4 | Scenario 5 |
|---|---|---|---|---|---|
| $t_0$ | 7.16 | 5.72 | 5.34 | 0 | 0 |
| $t_e$ | 16.36 | 14.7 | 14.36 | 8.52 | 8.42 |
| $t_e - t_0$ | 9.2 | 8.98 | 9.02 | 8.52 | 8.42 |

Table VI shows the start time $t_0$, the end time $t_e$ and the duration of the lane changing process in five scenarios. The start time $t_0$ is determined if the condition in Equation (11) is satisfied. While the end time of the lane-changing process $t_e$ can be obtained when the lateral tracking error reaches 0 as can be seen from Fig. 17.

The results of scenarios 1-3 show that reasonable parameter values can keep the lane-changing maneuver to be completed safely and efficiently. However, if the $TFV$'s speed changes, the trajectory planning module needs to re-plan a trajectory in line with the current state during the actual lane-changing process, otherwise the lane-changing could be unsafe. Fig. 18 (a) shows that the $SV$ begins to change lanes at $0s$, and the $TFV$ begins to decelerate at $0s$ in scenario 4. The detailed process of dynamic lane-changing trajectory planning of the $SV$ in scenario 4 is illustrated in Fig 18 (b). The proposed model updates the current optimal lane-changing trajectory at every sampling time, and the final lane-changing trajectory passes through the planned trajectory of all sampling times.

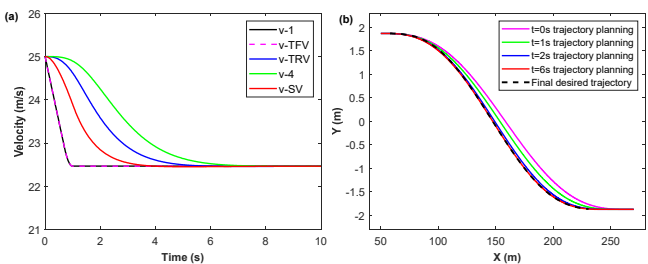

Fig. 18. The dynamic lane-changing trajectory planning process in scenario 4

## VI. CONCLUSION

In this paper, we propose a dynamic cooperative vehicle platooning approach considering both the longitudinal and the lateral control. The DCPID algorithm is developed to provide the longitudinal control for vehicle platooning. This algorithm can not only allow the platoon running steadily with slight disturbance but can also quickly return to the steady state under the extremely unstable condition. The strong stability and good anti-interference performance make the proposed algorithm very promising for real world implementation. Moreover, the simulation results show that the proposed DCPID approach is superior to the DMPC and the single PID methods in terms of anti-interference ability.

We further propose the cooperative lane-changing model by combining the DCPID algorithm with the sine function. The longitudinal acceleration and speed of the lane-changing vehicle are determined by the DCPID algorithm considering the speed variations of the vehicle in front on the target lane ($TFV$). The reference trajectory is planned using the sine function taking into account of ride comfortability at different speeds and can be updated in real time in response to the state change of the $TFV$, such as emergency operations (sudden deceleration). Several numerical simulations were carried out to demonstrate the performance of the proposed method. The results show that the lane-changing trajectories generated by our proposed cooperative lane-changing algorithm are rather smooth and highly consistent with the desired trajectories at different speeds.

This paper focuses on the control strategies of CAV platooning on the basic section of the freeway. In future, we would like to continue our research on developing control strategies of vehicle platooning in more complex situations including on-ramps and off-ramps, as well as under different traffic flow conditions such as the free flow and the transition from the free flow to the congested flow. In addition, we would also like to extend the proposed approach to the mixed traffic flow with CAVs and human-driven vehicles and investigate the impact of the CAV penetration rate on the performance of the mixed traffic flow.